\documentclass[12pt]{article}
\usepackage[pagebackref]{hyperref}

\usepackage[utf8]{inputenc}
\usepackage[T1]{fontenc}

\usepackage[]{amsfonts}
\usepackage{amssymb}
\usepackage{amsmath,amsthm}
\def\couleur(#1 #2 #3)
	{
\special{ps::[end]}}

\def\bx#1{\setbox1=\hbox{\kern3pt{#1}\kern3pt}			
 \dimen1=\ht1 \advance\dimen1 by 3pt \dimen2=\dp1 \advance\dimen2 by 3pt
 \setbox1=\hbox{\vrule height\dimen1 depth\dimen2\box1\vrule}%
 \setbox1=\vbox{\hrule\box1\hrule}%
 \advance\dimen1 by .4pt \ht1=\dimen1
 \advance\dimen2 by .4pt \dp1=\dimen2 \box1\relax}

\def\wbb#1{\kern#1em}
\def\vci{\vrule  width.02em height1.47ex depth-.0ex}		
\def\11{{\rm\wbb{.2}\vci\wbb{-.37}1}}

\def\underset#1#2{\mathrel{\mathop{\kern0pt #2}\limits_{#1}}}

\def\overset#1#2{\mathrel{\mathop{\kern0pt #2}\limits^{#1}}}

\newtheorem{thm}{Theorem}[section]
\newtheorem{lem}[thm]{Lemma}
\newtheorem{cor}[thm]{Corollary}

\begin{document}

\title{The raising steps method. Application to the $\displaystyle L^{r}$ Hodge theory in a compact riemannian manifold.}

\author{Eric Amar}

\date{}
\maketitle
 \renewcommand{\abstractname}{Abstract}

\begin{abstract}
Let $X$ be a complete metric space and $\displaystyle \Omega
 $ a domain in $\displaystyle X.$ The Raising Steps Method allows
 to get from local results on solutions $u$ of a linear equation
 $\displaystyle Du=\omega $ global ones in $\displaystyle \Omega .$\ \par 
It was introduced in~\cite{AmarSt13} to get good estimates on
 solutions of $\bar \partial $ equation in domains in a Stein manifold.\ \par 
As a simple application we shall get a strong $\displaystyle
 L^{r}$ Hodge decomposition theorem for $p-$forms in a compact
 riemannian manifold without boundary, and then we retrieve this
 known result by an entirely different and simpler method.\ \par 
\end{abstract}

\section{Introduction.}
\quad This work proposes a way for passing from local to global: the
 raising steps method, RSM for short. I introduce it precisely
 to get $L^{r}-L^{s}$ estimates for solutions of the $\bar \partial
 $ equation in Stein manifold in~\cite{AmarSt13}. See also~\cite{dBarIntersection17}
 to have result in case of intersection of domains in Stein manifolds.\ \par 
\quad The aim is to generalise it to the case where the $\bar \partial
 $ operator is replaced by an abstract linear operator $D$ acting
 on a domain in a complete metric space.\ \par 
\quad We shall deal with the following situation: we have a complete
 metric space $X$ admitting partitions of unity (see condition
 (ii) below) and a measure $\mu .$\ \par 
\quad A \emph{domain} of $X$ will be a connected open set $\Omega $
 of $X,$ relatively compact.\ \par 
\quad We are interested on solutions $u$ of a linear equation $Du=\omega
 ,$ in a domain $\Omega .$ Precisely fix a threshold $s>1.$ Suppose
 you have a \emph{global} solution $\displaystyle u$ on $\Omega
 $ of $Du=\omega $ with estimates $L^{s}(\Omega )\rightarrow
 L^{s}(\Omega ).$ It may happen that we have a constrain:\ \par 
\quad \quad \quad $\exists K\subset L^{s'}(\Omega )$ with $s'$ the conjugate exponent
 for $s,$ such that $\forall u::Du\in L^{s}(\Omega ),\ \forall
 h\in K,\ {\left\langle{Du,h}\right\rangle}=0.$\ \par 
In this case, in order for a solution of $Du=\omega $ to exist,
 we need to have $\omega \perp K.$ With no constrain, we take
 $K=\lbrace 0\rbrace .$\ \par 
\quad Very often this threshold will be $s=2,$ since Hilbert spaces
 are usually more tractable.\ \par 
\quad Now suppose that we have, for $1\leq r\leq s,$ \emph{local} solutions
 $\displaystyle u$ on $U\cap \Omega $ of $Du=\omega $ with estimates
 $L^{r}(\Omega )\rightarrow L^{t}(U\cap \Omega )$ with a \emph{strict
 increase} of the regularity, for instance $\displaystyle \frac{1}{t}=\frac{1}{r}-\tau
 ,\ \tau >0$ for any $\displaystyle r\leq s,$ then the raising
 steps method gives a \emph{global} solution $v$ of $Dv=\omega
 $ which is essentially in $L^{t}(\Omega )$ if we start with
 a data $\omega $ in $L^{r}(\Omega ).$\ \par 
\quad In particular we prove:\ \par 

\begin{thm}
(Raising steps theorem) Under the assumptions above, there is
 a positive constant $\displaystyle c$ such that for $1\leq r\leq
 s,$ if $\omega \in L^{r}(\Omega ),\ \omega \perp K,$ there is
 a$\ u\in L^{t}(\Omega )$ with $\displaystyle \frac{1}{t}=\frac{1}{r}-\tau
 ,$ such that $Du=\omega +\tilde \omega ,$ with $\tilde \omega
 \in L^{s}(\Omega ),\ \tilde \omega \perp K$ and control of the norms.
\end{thm}
\quad From this theorem and the fact that there is a \emph{global}
 solution $v\in L^{s}(\Omega )$ to $Dv=\bar \omega $ we get that
 $u-v$ is the \emph{global} solution we are searching for.\ \par 
\ \par 
\quad To illustrate the method, we shall apply it for the Poisson equation
 associated to the Hodge Laplacian in a compact riemannian boundary-less
 manifold $(M,g).$\ \par 
On $(M,g)$ we can define Sobolev spaces $W^{k,r}(\Omega )$ (see~\cite{Hebey96})
 and if $M$ is compact these spaces are in fact independent of
 the metric. Moreover the Sobolev embeddings are true in this
 case and a chart diffeomorphism makes a correspondence between
 Sobolev spaces in ${\mathbb{R}}^{n}$ and Sobolev spaces in $M.$\ \par 
\quad Let $d$ be the exterior derivative on $M$ and $d^{*}$ its adjoint;
 we define the Hodge laplacian acting from $p$  differential
 forms to $p$ differential forms to be: $\Delta :=dd^{*}+d^{*}d.$
 Because $M$ is compact, we have that $\Delta $ is self adjoint.
 The Poisson equation on $M$ is, for a given $p$-form $\omega
 $ on $M,$ to find a $p$-form $u$ on $M$ such that $\Delta u=\omega .$\ \par 
\quad Let ${\mathcal{H}}_{p}$ be the set of $p$-\emph{harmonic} forms
 in $M,$ i.e. $h\in {\mathcal{H}}_{p}\iff h\in {\mathcal{C}}^{\infty
 }_{p}(M),\ \Delta h=0.$ We have:\ \par 
\quad \quad \quad $\forall h\in {\mathcal{H}}_{p},\ {\left\langle{\Delta u,h}\right\rangle}={\left\langle{u,\Delta
 h}\right\rangle}=0$\ \par 
hence, in order to solve $\Delta u=\omega ,$ we need to have
 $\omega \perp {\mathcal{H}}_{p}.$\ \par 
\quad We derive, from a solution of the Poisson equation we get by
 use of the RSM, a $L^{r}_{p}$ Hodge decomposition for $p$ differential
 forms on $M.$\ \par 
\quad We shall use, for the local results, the classical ones. Let
 $B$ be a ball in ${\mathbb{R}}^{n},$ then:\ \par 
\quad \quad \quad $\displaystyle \forall \gamma \in L^{r}(B),\ \exists u\in W^{2,r}(B)::\Delta
 _{{\mathbb{R}}}u=\gamma ,\ {\left\Vert{v_{0}}\right\Vert}_{W^{2,r}(B)}\leq
 C{\left\Vert{\gamma }\right\Vert}_{L^{r}(B)}.$\ \par 
These non trivial estimates are coming from Gilbarg and Trudinger~\cite[Theorem
 9.9, p. 230]{GuilbargTrudinger98} and the constant $C=C(n,r)$
 depends only on $n$ and $r.$\ \par 
\quad Then, together with the R.S.M., we get a solution of the Poisson
 equation for the Hodge laplacian.\ \par 

\begin{thm}
Let $M$ be a compact, ${\mathcal{C}}^{\infty }$ Riemannian manifold
 without boundary. For any $r,\ 1\leq r\leq n/2,$ if $g$ is a
 $p$-form in $L_{p}^{r}(M)\cap {\mathcal{H}}_{p}^{\perp }$ there
 is a $p$-form $v\in L^{t}(M)$ such that $\Delta v=g$ and ${\left\Vert{v}\right\Vert}_{L^{t}_{p}(M)}\leq
 c{\left\Vert{g}\right\Vert}_{L^{r}_{p}(M)}$ with $\displaystyle
 \frac{1}{t}=\frac{1}{r}-\frac{2}{n}.$\par 
Moreover for any $r>1$ and any $p$ form $u$ solution of $\Delta
 u=g,$ we have $u$ in $W_{p}^{2,r}(M).$
\end{thm}
\quad From this result we deduce the $L^{r}$ Hodge decomposition of $p$-forms.\ \par 

\begin{thm}
Let $(M,g)$ be a compact riemannian manifold without boundary.
 We have the strong $L^{r}$ Hodge decomposition:\par 
\quad \quad $\displaystyle \forall r,\ 1\leq r<\infty ,\ L^{r}_{p}(M)={\mathcal{H}}_{p}^{r}\oplus
 \mathrm{I}\mathrm{m}\Delta (W^{2,r}_{p}(M))={\mathcal{H}}_{p}^{r}\oplus
 \mathrm{I}\mathrm{m}d(W^{1,r}_{p}(M))\oplus \mathrm{I}\mathrm{m}d^{*}(W^{1,r}_{p}(M)).$
\end{thm}

      The case $r=2$ of this decomposition which gives us $s=2$
 as a threshold, goes to Morrey in 1966 and essentially all results
 in the $L^{2}$ case we use here are coming from the basic work
 of Morrey~\cite{Morrey66}.\ \par 
\quad This decomposition is an already known result of C. Scott~\cite{Scott95}
 but proved here by an entirely different method. Critical to
 Scott's proof is a nice $L^{r}$ Gaffney's inequality which he
 proved and used to get the $L^{r}$ Hodge decomposition, the
 same way than Morrey~\cite{Morrey66} did with the $L^{2}$ Gaffney's
 inequality~\cite{Gaffney55} to get the $L^{2}$ Hodge decomposition.\ \par 
\quad In the case of a compact manifold with boundary, G. Schwarz~\cite{Schwarz95}
 proved also a $L^{r}$  Gaffney's inequality to get the $L^{r}$
 Hodge decomposition in that case, then he deduced of it a global
 $L^{r}$ solution for the equation $\Delta u=\omega .$\ \par 
\quad In the nice book by F.W. Warner~\cite{Warner83}, the author proved
 directly, without the use of Gaffney's inequality, a global
 $L^{2}$ solution for the equation $\Delta u=\omega $ in the
 case of a compact manifold without boundary. He deduced from
 it the $L^{2}$ Hodge decomposition.\ \par 
\quad Here we use the RSM plus the global  $L^{2}$ solution for the
 equation $\Delta u=\omega $ given by Warner~\cite{Warner83},
 to get a global  $L^{r}$ solution for the equation $\Delta u=\omega
 $ and then recover the $L^{r}$ Hodge decomposition, without
 any Gaffney's inequalities. Hence we get a completely different
 proof of the known  $L^{r}$ Hodge decomposition.\ \par 
\quad Many important applications of the $L^{2}$ Hodge decomposition
 in cohomology theory and algebraic geometry are in the book
 by C. Voisin~\cite{Voisin02}.\ \par 
\quad So it may be interesting to have a short proof of this important
 Hodge decomposition in the $L^{r}$ case.\ \par 
\quad Finally, in the last section we prove, by use of the "double
 manifold" technique:\ \par 

\begin{thm}
Let $\Omega $ be a domain in the smooth complete riemannian manifold
 $M$ and $\omega \in L^{r}_{p}(\Omega ),$ then there is a $p$-form
 $u\in W^{2,r}_{p}(\Omega ),$ such that $\Delta u=\omega $ and
 ${\left\Vert{u}\right\Vert}_{W^{2,r}_{p}(\Omega )}\leq c(\Omega
 ){\left\Vert{\omega }\right\Vert}_{L^{r}_{p}(\Omega )}.$
\end{thm}
\quad And we make a short incursion in the domain of manifold with boundary:\ \par 

\begin{cor}
Let $M$ be a smooth compact riemannian manifold with smooth boundary
 $\partial M.$ Let $\displaystyle \omega \in L^{r}_{p}(M).$ 
 There is a $p$-form $u\in W^{2,r}_{p}(M),$ such that $\Delta
 u=\omega $ and ${\left\Vert{u}\right\Vert}_{W^{2,r}_{p}(M)}\leq
 c{\left\Vert{\omega }\right\Vert}_{L^{r}_{p}(M)}.$
\end{cor}
\quad Schwarz~\cite[Theorem 3.4.10, p. 137]{Schwarz95} proved a better
 theorem: you can prescribe the values of $u$ on the boundary.
 But again the proof here is much lighter.\ \par 

\section{The Raising Steps Method.}
\quad \quad 	We shall deal with the following situation: we have a complete
 metric space $X$ admitting partitions of unity (see condition
 (ii) below) and a positive $\sigma $-finite measure $\mu .$\ \par 

\subsection{Assumptions on the linear operator $D.$}
\quad We shall denote $E^{p}(X)$ the set of ${\mathbb{C}}^{p}$ valued
 fonctions on $X.$ This means that $\omega \in E^{p}(X)\iff \omega
 (x)=(\omega _{1}(x),...,\omega _{p}(x)).$ We put a punctual
 norm on $\omega \in E^{p}(X),\ \left\vert{\omega (x)}\right\vert
 ^{2}:=\sum_{j=1}^{p}{\left\vert{\omega _{j}(x)}\right\vert ^{2}}$
 and if $\displaystyle U$ is an open set in $X,$ we consider
 the Lebesgue space $L^{r}_{p}(U)$:\ \par 
\quad \quad \quad $\displaystyle \omega \in L^{r}_{p}(U)\iff {\left\Vert{\omega
 }\right\Vert}^{r}_{L^{r}_{p}(U)}:=\int_{U}{\left\vert{\omega
 (x)}\right\vert ^{r}d\mu (x)}<\infty .$\ \par 
The space $\displaystyle L^{2}_{p}(U)$ is a Hilbert space with
 the scalar product ${\left\langle{\omega ,\omega '}\right\rangle}:=\int_{U}{{\left({\sum_{j=1}^{p}{\omega
 _{j}(x)\bar \omega '_{j}(x)}}\right)}d\mu (x)}.$\ \par 
\quad We are interested in solution of a linear equation $Du=\omega
 ,$ where $D=D_{p}$ is a linear operator acting on $E^{p}.$\ \par 
\quad In order to have $\displaystyle Du=\omega ,$ it may happen that
 we have a constrain: there is a subspace $K\subset L^{r'}(\Omega
 )$ such that\ \par 
\quad \quad \quad $\forall h\in K,\ \forall u::Du\in L^{r}(\Omega ),\ {\left\langle{Du,h}\right\rangle}=0\iff
 Du\perp K.$\ \par 
The absence of constrain is done by setting $K=\lbrace 0\rbrace .$\ \par 
\quad Now on the integer $p$ will be fixed so the explicit mention
 of the integer $p$ will be often omitted.\ \par 
\quad \quad 	We shall make the following hypotheses.\ \par 
\quad \quad 	Let $\Omega $ be a domain in $X.$ There is a $\tau \geq \delta
 $ with $\displaystyle \frac{1}{t}=\frac{1}{r}-\tau ,$ and a
 positive constant $c_{l}$ such that:\ \par 
\quad \quad \quad (i) Local Existence with Increasing Regularity (LEIR):  for any
 $x\in \bar \Omega ,$ there is a ball $B:=B(x,R_{x})$ such that
 if $\omega \in L_{p}^{r}(B),$ we can solve $Du_{x}=\omega $
 in $B':=B(x,R_{x}/2)$ with $L_{p}^{r}(B)-L_{p}^{t}(B')$ estimates,
 i.e.  	$\exists u_{x}\in L^{t}(B'),\ Du_{x}=\omega $ in $B'$
 and ${\left\Vert{u_{x}}\right\Vert}_{L^{t}(B')}\leq c_{l}{\left\Vert{\omega
 }\right\Vert}_{L^{r}(B)}.$\ \par 
\quad It may append, in the case $X$ is a manifold, that we have a
 better regularity for the local existence:\ \par 
\quad \quad \quad (i') Sobolev regularity: if $\omega \in L_{p}^{r}(B),$ we can
 solve $Du_{x}=\omega $ in $B':=B(x,R_{x}/2)$ with $L_{p}^{r}(B)-W_{p}^{\alpha
 ,r}(B')$ estimates, i.e.  	$\exists u_{x}\in W_{p}^{\alpha ,r}(B'),\
 Du_{x}=\omega $ in $B'$ and ${\left\Vert{u_{x}}\right\Vert}_{W_{p}^{\alpha
 ,r}(B')}\leq c_{l}{\left\Vert{\omega }\right\Vert}_{L^{r}(B)}.$\ \par 
\quad By compactness we can cover $\bar \Omega $ by a finite set of
 balls $\displaystyle \lbrace B(x_{j},R_{j}/2)\rbrace _{j=1,...,N}$
 of the previous form. Set $B_{j}:=B(x_{j},R_{j}),\ B'_{j}:=B(x_{j},R_{j}/2).\
 $Set $u_{j}$ the local solution of $\displaystyle Du_{j}=\omega
 $ with $\displaystyle {\left\Vert{u_{j}}\right\Vert}_{L^{t}(B_{j}')}\leq
 c_{l}{\left\Vert{\omega }\right\Vert}_{L^{r}(B_{j})}.$\ \par 
\quad \quad \quad (ii) Partition of unity: If $\lbrace B'_{j}\rbrace _{j=1,...,N}$
 is a covering of $\bar \Omega ,$ then there is an associated
 set of functions  $\displaystyle \lbrace \chi _{j}\rbrace _{j=1,...,N}$
 such that $\chi _{j}$ has compact support in $B'_{j},\ \forall
 j=1,...,N,\ 0\leq \chi _{j}(x)\leq 1$ and $\sum_{j=1}^{N}{\chi
 _{j}(x)}=1$ for $x\in \bar \Omega .$\ \par 
\quad \quad \quad (iii) Commutator condition: We set $\Delta _{j}=\Delta (\chi
 _{j},u_{j}):=\chi _{j}Du_{j}-D(\chi _{j}u_{j}).$ There is a
 constant $\delta >0$ such that, with $\displaystyle \frac{1}{t}=\frac{1}{r}-\delta
 ,$ we have:\ \par 
\quad \quad \quad \quad \quad $\displaystyle {\left\Vert{\Delta _{j}}\right\Vert}_{L^{t}(B_{j}')}\leq
 c(\chi _{j})({\left\Vert{\omega }\right\Vert}_{L^{r}(B_{j})}+{\left\Vert{u_{j}}\right\Vert}_{L^{r}(B_{j})}).$\
 \par 
\quad \quad \quad  (iv) Global resolvability: We can solve $Dw=\omega $ globally
 in $\Omega $ with $L^{s}-L^{s}$ estimates, i.e.\ \par 
\quad \quad \quad \quad 	$\exists c_{g}>0,\ \exists w\ s.t.\ Dw=\omega $ in $\displaystyle
 \Omega $ and $\ {\left\Vert{w}\right\Vert}_{L^{s}(\Omega )}\leq
 c_{g}{\left\Vert{\omega }\right\Vert}_{L^{s}(\Omega )},$ provided
 that $\omega \perp K.$\ \par 
\quad It may append, in the case $X$ is a manifold, that we have a
 better regularity for the global existence:\ \par 
\quad \quad (iv') Sobolev regularity: We can solve $Dw=\omega $ globally
 in $\Omega $ with $L_{p}^{s}-W_{p}^{\alpha ,s}$ estimates, i.e.\ \par 
\quad \quad \quad \quad 	$\exists c_{g}>0,\ \exists w\ s.t.\ Dw=\omega $ in $\displaystyle
 \Omega $ and $\ {\left\Vert{w}\right\Vert}_{W_{p}^{\alpha ,s}(\Omega
 )}\leq c_{g}{\left\Vert{\omega }\right\Vert}_{L^{s}(\Omega )},$
 provided that $\displaystyle \omega \perp K.$\ \par 

\begin{thm}
~\label{CC2}(Raising steps theorem) Under the assumptions (i),
 (ii), (iii), (iv) above, there is a positive constant $\displaystyle
 c_{f}$ such that for $1\leq r\leq s,$ if $\omega \in L^{r}(\Omega
 ),\ \omega \perp K$ there is a$\ u=u_{s}\in L^{t}(\Omega )$
 with $\displaystyle \frac{1}{t}=\frac{1}{r}-\tau ,$ such that
 $Du=\omega +\tilde \omega ,$ with $\tilde \omega \in L^{s}(\Omega
 ),\ \tilde \omega \perp K$ and control of the norms.\par 
\quad If moreover we have (i') then $u\in W_{p}^{\alpha ,r}(\Omega
 )$ with control of the norm.
\end{thm}

\begin{proof}
Let $r\leq s$ and $\omega \in L^{r}(\Omega ),\ \omega \perp K;$
 we start with the covering $\lbrace B'_{j}\rbrace _{j=1,...,N}$
 and the local solution $Du_{j}=\omega $ with $\displaystyle
 \frac{1}{t}=\frac{1}{r}-\tau ,\ u_{j}\in L^{t}(B'_{j})$ given
 by hypothesis (i).\par 
If (i') is true, then we have $\displaystyle u_{j}\in W_{p}^{\alpha
 ,r}(B_{j}'),\ Du_{j}=\omega .$\par 
Let $\lbrace \chi _{j}\rbrace _{j=1,...,N}$ be the partition
 of unity subordinate to $\lbrace B'_{j}\rbrace _{j=1,...,N}$
 given by (ii). Because $\displaystyle 0\leq \chi _{j}(x)\leq
 1$ we have $\chi _{j}u_{j}\in L^{t}(\Omega )$ hence\par 
\quad \quad \quad \quad \quad  \begin{equation} {\left\Vert{ \chi _{j}u_{j}}\right\Vert}_{L^{t}(\Omega
 )}\leq {\left\Vert{u_{j}}\right\Vert}_{L^{t}(B'_{j})}\leq c_{l}{\left\Vert{\omega
 }\right\Vert}_{L^{r}(\Omega )}.\label{rS25}\end{equation}\par 
Set 	$v_{0}:=\sum_{j=1}^{N}{\chi _{j}u_{j}}.$ Then we have, setting
 now $\displaystyle \frac{1}{t_{0}}=\frac{1}{r}-\tau $:\par 
\quad \quad 	$\bullet $ $v_{0}\in L^{t_{0}}(\Omega )$ because  $\chi _{j}u_{j}\in
 L^{t_{0}}(\Omega )$ for $j=1,...,N,$ and  ${\left\Vert{v_{0}}\right\Vert}_{L^{t_{0}}(\Omega
 )}\leq c{\left\Vert{\omega }\right\Vert}_{L^{r}(\Omega )}$ with
  $c=Nc_{l}$ by~(\ref{rS25}). If (i') is true, i.e. $X$ is a
 manifold, then we can choose $\chi _{j}\in {\mathcal{D}}(B'_{j}),$
 i.e. in the space of ${\mathcal{C}}^{\infty }$ functions of
 compact support in $B_{j},$ hence ${\left\Vert{\chi _{j}u_{j}}\right\Vert}_{W_{p}^{\alpha
 ,r}(\Omega )}\leq c{\left\Vert{u_{j}}\right\Vert}_{W_{p}^{\alpha
 ,r}(B'_{j})}\leq c_{l}{\left\Vert{\omega }\right\Vert}_{L^{r}(\Omega
 )}.$ So\par 
\quad \quad $\bullet $ $\displaystyle v_{0}\in W_{p}^{\alpha ,r}(\Omega )$
 with $\displaystyle {\left\Vert{v_{0}}\right\Vert}_{W_{p}^{\alpha
 ,r}(\Omega )}\leq c{\left\Vert{\omega }\right\Vert}_{L^{r}(\Omega )}.$\par 
We have\par 
\quad \quad 	$\bullet $ $Dv_{0}=\sum_{j=1}^{N}{\chi _{j}Du_{j}}+\sum_{j=1}^{N}{\Delta
 (\chi _{j},u_{j})}.$\par 
Setting  $\displaystyle \omega _{1}(x):=\sum_{j=1}^{N}{\Delta
 (\chi _{j},u_{j})(x)},$ we get\par 
\quad \quad \quad $\displaystyle \forall x\in \Omega ,\ Dv_{0}(x)=\sum_{j=1}^{N}{\chi
 _{j}(x)\omega (x)}+\sum_{j=1}^{N}{\Delta (\chi _{j},u_{j})(x)}=\omega
 (x)+\omega _{1}(x)$\par 
\quad By hypothesis (iii), $\Delta (\chi _{j},u_{j})\in L^{s_{0}}(B'_{j})$
 with $\displaystyle \frac{1}{s_{0}}=\frac{1}{r}-\delta ,$ hence
 $\omega _{1}\in L^{s_{0}}(\Omega ),$ with\par 
\quad \quad \quad \begin{equation} {\left\Vert{ \omega _{1}}\right\Vert}_{L^{s_{0}}(\Omega
 )}\leq G{\left\Vert{\omega }\right\Vert}_{L^{r}(\Omega )}\label{lR3}\end{equation}\par
 
and $G=c_{sl}N.$\par 
\quad The regularity of $\omega _{1}$ is higher by one step $\delta
 >0$ than that of $\omega .$ Moreover\par 
\quad \quad \quad $\forall h\in K,\ {\left\langle{\omega _{1},h}\right\rangle}={\left\langle{\omega
 ,h}\right\rangle}-{\left\langle{Dv_{0},h}\right\rangle}=0$ because
 $\displaystyle {\left\langle{\omega ,h}\right\rangle}=0$ and
 $\displaystyle {\left\langle{Dv_{0},h}\right\rangle}=0.$\par 
\quad 	If $s_{0}\geq s$ we notice that $\omega _{1}\in L^{s}(\Omega
 )$ because $L^{s_{0}}(\Omega )\subset L^{s}(\Omega )$ for $\Omega
 $ is relatively compact. So we are done by setting $u_{s}=v_{0},\
 \tilde \omega =\omega _{1}.$\par 
\quad If $s_{0}<s$ we proceed by induction: we set $t_{1}$ such that
 $\displaystyle \frac{1}{t_{1}}=\frac{1}{s_{0}}-\tau =\frac{1}{r}-\delta
 -\tau ,$ and we use the same covering  $\lbrace B'_{j}\rbrace
 _{j=1,...,N}$ and the same partition of unity $\lbrace \chi
 _{j}\rbrace _{j=1,...,N}$ and with $\omega _{1}$ in place of
 $\omega ,$ we get 		$\exists v_{1}\in L^{t_{1}}(\Omega ),\ Dv_{1}=\omega
 _{1}+\omega _{2}$ and, setting $s_{1}$ such that $\displaystyle
 \frac{1}{s_{1}}=\frac{1}{s_{0}}-\delta =\frac{1}{r}-2\delta
 ,$ we have that the regularity of $\omega _{2}$ raises of $2$
 times $\delta $ from that of $\omega .$\par 
\quad We still have\par 
\quad \quad \quad $\displaystyle \forall h\in K,\ {\left\langle{\omega _{2},h}\right\rangle}={\left\langle{\omega
 _{1},h}\right\rangle}-{\left\langle{Dv_{1},h}\right\rangle}=0$\par 
so by induction, after a finite number $k$ of steps, we get a
 $s_{k}\geq s.$ The linear combination $u:=\sum_{j=0}^{k-1}{(-1)^{j}v_{j}}$
 now gives $Du=\omega +(-1)^{k}\omega _{k}$ with $\omega _{k}\perp
 K.$ And again we are done by setting $\tilde \omega =(-1)^{k}\omega _{k}.$\par 
\quad If (i') is true, then we have $v_{0}\in W^{\alpha ,r}(\Omega
 )$ and $v_{1}\in W^{\alpha ,s_{0}}(\Omega ),\ {\left\Vert{v_{1}}\right\Vert}_{W^{\alpha
 ,s_{0}}(\Omega )}\leq c{\left\Vert{\omega _{1}}\right\Vert}_{L^{s_{0}}(\Omega
 )}.$ Hence by~(\ref{lR3}) we get ${\left\Vert{v_{1}}\right\Vert}_{W^{\alpha
 ,s_{0}}(\Omega )}\leq cG{\left\Vert{\omega }\right\Vert}_{L^{r}(\Omega
 )}.$\par 
A fortiori, $v_{1}\in W^{\alpha ,r}(\Omega )$ because $r<s_{0}$
 with the same control of the norm. Likewise we have $\forall
 j=1,...,k,\ {\left\Vert{v_{j}}\right\Vert}_{W^{\alpha ,r}(\Omega
 )}\leq c{\left\Vert{\omega }\right\Vert}_{L^{r}(\Omega )},$
 hence the same control of the norm of $u.$
\end{proof}

\begin{cor}
~\label{CC3}Under the assumptions of the raising steps theorem
 and with the global assumption (iv), there is a constant $c_{f}>0,$
 such that for $r\leq s,$ if $\omega \in L^{r}(\Omega ),\ \omega
 \perp K$ there is a$\ v\in L^{t}(\Omega )$ with $\displaystyle
 t:=\min (s,t_{0})$ and $\displaystyle \frac{1}{t_{0}}=\frac{1}{r}-\tau
 ,$ such that\par 
\quad \quad \quad \quad 	$Dv=\omega $ and $v\in L^{t}(\Omega ),\ {\left\Vert{v}\right\Vert}_{L^{t}(\Omega
 )}\leq c{\left\Vert{\omega }\right\Vert}_{L^{r}(\Omega )}.$\par 
If moreover (i') and (iv') are true then we have $\displaystyle
 v\in W^{\alpha ,r}(\Omega )\cap L^{t}(\Omega )$ with $\displaystyle
 {\left\Vert{v}\right\Vert}_{W^{\alpha ,r}(\Omega )}\leq c{\left\Vert{\omega
 }\right\Vert}_{L^{r}(\Omega )}.$
\end{cor}

\begin{proof}
By the raising steps Theorem~\ref{CC2} we have a $u\in L^{t_{0}}(\Omega
 )$ such that $Du=\omega +\tilde \omega $ with $\tilde \omega
 \in L^{s}(\Omega )$ and $\tilde \omega \perp K;$ by hypothesis
 (iv) we can solve $D\tilde v=\tilde \omega $ with $\tilde v\in
 L^{s}(\Omega )$ so it remains to set $v:=u-\tilde v$ to get
 $v\in L^{t}(\Omega )$ and $Dv=\omega .$\par 
\quad If moreover (i') is true then we have $\displaystyle u\in W^{\alpha
 ,r}(\Omega ).$ If (iv') is true then we can solve $\displaystyle
 D\tilde v=\tilde \omega $ with $\displaystyle \tilde v\in W^{\alpha
 ,s}(\Omega ).$ Hence, because $r\leq s,$ with $v:=u-\tilde v$
 we get $v\in W^{\alpha ,r}(\Omega )$ and $Dv=\omega .$ The proof is complete.
\end{proof}

\section{Application to Poisson equation on a compact riemannian manifold.}

\subsection{Local existence with increasing regularity.}
\quad In order to have the local result, we choose a chart $(V,\ \varphi
 :=(x_{1},...,x_{n}))$ so that $g_{ij}(y)=\delta _{ij}$ and $\varphi
 (V)=B$ where $B=B_{e}$ is an Euclidean ball centered at $\varphi
 (y)=0$ and $g_{ij}$ are the components of the metric tensor
 w.r.t. $\varphi .$\ \par 
\quad Because I was unable to find an easy proof of the following theorem
 in the literature, I reprove it for the reader's convenience.\ \par 

\begin{thm}
~\label{RSMA0} For any $y\in M,$ there are open sets  $W\subset
 \bar W\subset V\subset M,\ y\in W,$ such that we have:\par 
\quad \quad \quad $\displaystyle \forall \omega \in L_{p}^{r}(W),\ \exists u\in
 W_{p}^{2,r}(W)::\Delta u=\omega ,\ {\left\Vert{u}\right\Vert}_{W^{2,r}(W)}\leq
 C{\left\Vert{\omega }\right\Vert}_{L^{r}(W)}.$
\end{thm}

\begin{proof}
Of course the operator $d$ on $p$-forms is local and so is $d^{*}$
 as a first order differential operator.\par 
We start with a chart $(V,\varphi )$ of $M$ such that $\varphi
 (y)=0\in {\mathbb{R}}^{n}$ and the  metric tensor read in this
 chart at $y$ is the identity.\par 
\quad Then the Hodge laplacian $\Delta _{\varphi }$ read by $\varphi
 $ in a ball $B:=B(0,R)\subset {\mathbb{R}}^{n}$ is not so different
 from that of ${\mathbb{R}}^{n}$ in $B$ when acting on $p$-forms
 in $B.$ We set $\Delta _{\varphi }\omega _{\varphi }=\Delta
 _{{\mathbb{R}}}\omega _{\varphi }+A\omega _{\varphi },$ where
 $\omega _{\varphi }$ is the $p$-form $\omega $ read in the chart
 $(V,\varphi )$ and $A$ is a matrix valued second order operator
 with ${\mathcal{C}}^{\infty }$ smooth coefficients such that
 $A:W^{2,r}(B)\rightarrow L^{r}(B)$ with, for a $R$ small enough
 ${\left\Vert{Av}\right\Vert}_{L^{r}(B)}\leq c{\left\Vert{v}\right\Vert}_{W^{2,r}(B)}.$\par
 
\quad This is true because at the point $y\in V$ we are in the flat
 case and if $R$ is small enough, the difference $A:=\Delta _{\varphi
 }-\Delta _{{\mathbb{R}}}$ in operator norm $W^{2,r}(B)\rightarrow
 L^{r}(B)$ goes to $0$ when $R$ goes to $\displaystyle 0,$ because
 $\varphi \in {\mathcal{C}}^{2}$ and the metric tensor $g$ is
 also ${\mathcal{C}}^{2}.$\par 
\quad We know that $\Delta _{{\mathbb{R}}}$ operates component-wise
 on the $p$-form $\gamma ,$ so we have\par 
\quad \quad \quad $\forall \gamma \in L^{r}_{p}(B),\ \exists v_{0}\in W^{2,r}_{p}(B)::\Delta
 _{{\mathbb{R}}}v_{0}=\gamma ,\ {\left\Vert{v_{0}}\right\Vert}_{W^{2,r}(B)}\leq
 C{\left\Vert{\gamma }\right\Vert}_{L^{r}(B)},$\par 
simply setting the component of $v_{0}$ to be the Newtonian potential
 of the corresponding component of $\gamma $ in $U.$ This way
 $v_{0}$ is linear with respect to $\gamma .$ These non trivial
 estimates are coming from Gilbarg and Trudinger~\cite[Thorem
 9.9, p. 230]{GuilbargTrudinger98} and the constant $C=C(n,r)$
 depends only on $n$ and $r.$\par 
\quad So we get $\Delta _{{\mathbb{R}}}v_{0}+Av_{0}=\gamma +\gamma _{1},$ with\par 
\quad \quad \quad $\gamma _{1}=Av_{0}\Rightarrow {\left\Vert{\gamma _{1}}\right\Vert}_{L^{r}(B)}\leq
 c{\left\Vert{v_{0}}\right\Vert}_{W^{2,r}(B)}\leq cC{\left\Vert{\gamma
 }\right\Vert}_{L^{r}(B)}.$\par 
\quad We solve again\par 
\quad \quad \quad $\displaystyle \exists v_{1}\in W^{2,r}_{p}(B)::\Delta _{{\mathbb{R}}}v_{1}=\gamma
 _{1},\ {\left\Vert{v_{1}}\right\Vert}_{W^{2,r}(B)}\leq C{\left\Vert{\gamma
 _{1}}\right\Vert}_{L^{r}(B)}=C^{2}c{\left\Vert{\gamma }\right\Vert}_{L^{r}(B)},$\par
 
and we set\par 
\quad \quad \quad $\displaystyle \gamma _{2}:=Av_{1}\Rightarrow {\left\Vert{\gamma
 _{2}}\right\Vert}_{L^{r}(B)}\leq c{\left\Vert{v_{1}}\right\Vert}_{W^{2,r}(B)}\leq
 C{\left\Vert{\gamma _{1}}\right\Vert}_{L^{r}(B)}\leq C^{2}c^{2}{\left\Vert{\gamma
 }\right\Vert}_{L^{r}(B)}.$\par 
\quad And by induction:\par 
\quad \quad \quad $\displaystyle \forall k\in {\mathbb{N}},\ \gamma _{k}:=Av_{k-1}\Rightarrow
 {\left\Vert{\gamma _{k}}\right\Vert}_{L^{r}(B)}\leq c{\left\Vert{v_{k-1}}\right\Vert}_{W^{2,r}(B)}\leq
 C{\left\Vert{\gamma _{k-1}}\right\Vert}_{L^{r}(B)}\leq C^{k}c^{k}{\left\Vert{\gamma
 }\right\Vert}_{L^{r}(B)}$\par 
and\par 
\quad \quad \quad $\displaystyle \exists v_{k}\in W^{2,r}_{p}(B)::\Delta _{{\mathbb{R}}}v_{k}=\gamma
 _{k},\ {\left\Vert{v_{k}}\right\Vert}_{W^{2,r}(B)}\leq C{\left\Vert{\gamma
 _{k}}\right\Vert}_{L^{r}(B)}\leq C^{k+1}c^{k}{\left\Vert{\gamma
 }\right\Vert}_{L^{r}(B)}.$\par 
Now we set $v:=\sum_{j\in {\mathbb{N}}}{(-1)^{j}v_{j}}.$ This
 series converges in norm $W^{2,r}(B),$ provided that we choose
 the radius of the ball $B$ small enough to have $cC^{2}<1,$ and we get:\par 
\quad \quad \quad $\displaystyle \Delta _{\varphi }v=\Delta _{{\mathbb{R}}}v+Av=\sum_{j\in
 {\mathbb{N}}}{(-1)^{j}(\Delta _{{\mathbb{R}}}v_{j}+Av_{j})}=\gamma ,$\par 
the last series converging in $L^{r}(B).$\par 
\quad In fact every step is linear and we get that $v$ is linear in $\gamma .$\par 
\quad Going back to the manifold $M$ with $\gamma :=\omega _{\varphi
 }$ and setting $u_{\varphi }:=v,\ W:=\varphi ^{-1}(B),$ we get
 the right estimates:\par 
\quad \quad \quad $\displaystyle \exists u\in W^{2,r}(W)::\Delta u=\omega \ in\
 W,\ {\left\Vert{u}\right\Vert}_{W^{2,r}(W)}\leq C{\left\Vert{\omega
 }\right\Vert}_{L^{r}(W)},$\par 
because the Sobolev spaces for $B$ go to the analogous Sobolev
 spaces for $W$ in $M.$
\end{proof}
\quad Now the next corollary is precisely the result we are searching for.\ \par 

\begin{cor}
For any $\displaystyle y\in M,$ there are open sets $V,\ W\Subset
 V,\ y\in W,$ such that we have for $r\leq 2$:\par 
$\displaystyle \forall \omega \in L^{r}_{p}(V),\ \exists u\in
 W^{2,r}(W)\cap L^{t}_{p}(W)::\Delta u=\omega ,\ {\left\Vert{u}\right\Vert}_{W^{2,r}(W)}\leq
 C{\left\Vert{\omega }\right\Vert}_{L^{r}(V)},\ {\left\Vert{u}\right\Vert}_{L^{t}_{p}(W)}\leq
 C{\left\Vert{\omega }\right\Vert}_{L^{r}_{p}(V)}$\par 
with $\displaystyle \frac{1}{t}=\frac{1}{r}-\frac{2}{n},$ and
 $\nabla u\in W^{1,r}_{p}(W).$
\end{cor}

\begin{proof}
The Theorem~\ref{RSMA0} gives $u\in W^{2,r}(W)$ such that $\Delta
 u=\omega ,\ {\left\Vert{u}\right\Vert}_{W^{2,r}(W)}\leq C{\left\Vert{\omega
 }\right\Vert}_{L^{r}(V)},$ hence we get that $\nabla u\in W^{1,r}_{p}(W)$
 with the same control: ${\left\Vert{\nabla u}\right\Vert}_{W^{1,r}(W)}\leq
 C{\left\Vert{\omega }\right\Vert}_{L^{r}(V)}.$\par 
\quad For the first statement it remains to apply the Sobolev embedding
 theorems which are true here.
\end{proof}
\quad So we are in a special case of the previous section with $D:=\Delta
 $ and, because $\Delta $ is essentially self adjoint, we have
 here $K:={\mathcal{H}}_{p}(M)$ where ${\mathcal{H}}_{p}(M)$
 is the space of ${\mathcal{C}}^{\infty }$ harmonic $p$-forms in $M.$\ \par 
\quad We shall need the following lemma.\ \par 

\begin{lem}
~\label{RSMH4}Let $\Delta _{\varphi }$ be a second order elliptic
 matrix operator with ${\mathcal{C}}^{\infty }$ coefficients
 operating on  $p$-forms $v$ defined in $U\subset {\mathbb{R}}^{n}.$
 Let $B:=B(0,R)$ a ball in ${\mathbb{R}}^{n},\ B':=B(0,R/2)$
 and suppose that $B\Subset U.$ Then we have an interior estimate:
 there are constants $c_{1},c_{2}$ depending only on $n=\mathrm{d}\mathrm{i}\mathrm{m}_{{\mathbb{R}}}M,\
 r$ and the ${\mathcal{C}}^{1}$norm of the coefficients of $\Delta
 _{\varphi }$ in $\bar B$ such that\par 
\begin{equation}  \forall v\in W_{p}^{2,r}(B),\ {\left\Vert{v}\right\Vert}_{W^{2,r}(B')}\leq
 c_{1}{\left\Vert{v}\right\Vert}_{L^{r}(B)}+c_{2}{\left\Vert{\Delta
 _{\varphi }v}\right\Vert}_{L^{r}(B)}.\label{CM11}\end{equation}
\end{lem}

\begin{proof}
For a $0$-form this lemma is exactly~\cite[Theorem 9.11]{GuilbargTrudinger98}.\par
 
\quad For $p$-forms we cannot avoid the use of deep results on elliptic
 systems of equations.\par 
\quad Let $\displaystyle v$ be a $p$-form in $B\subset {\mathbb{R}}^{n}.$
 We use the interior estimates in~\cite[\S  6.2, Thm 6.2.6]{Morrey66}.
 In our context, second-order elliptic system, and with our notations,
 with $r>1,$ we get:\par 
\[ \exists C>0,\ \forall v\in W_{p}^{2,r}(B),\ {\left\Vert{v}\right\Vert}_{W^{2,r}(B')}\leq
 c_{1}R^{-2}{\left\Vert{v}\right\Vert}_{L^{r}(B)}+c_{2}{\left\Vert{\Delta
 _{\varphi }v}\right\Vert}_{L^{r}(B)},\]\par 
already including the dependence in $R.$\par 
\quad The constants $c_{1},c_{2}$ depend only on $r,\ n:=\mathrm{d}\mathrm{i}\mathrm{m}M$
 and the bounds and moduli of continuity of all the coefficients
 of the matrix $\Delta _{\varphi }.$ (In~\cite{Morrey66}, p.
 213: the constant depends only on $E$ and on $E'.$)\par 
\quad In particular, if $\Delta _{\varphi }$ has its coefficients near
 those of $\Delta _{{\mathbb{R}}}$ in the ${\mathcal{C}}^{1}$
 norm, then the constants $c_{1},c_{2}$ are near the ones obtained
 for $\Delta _{{\mathbb{R}}}.$
\end{proof}
\quad Now we deduce from it local interior regularity for the laplacian
 on a smooth compact manifold without boundary.\ \par 

\begin{lem}
~\label{RSMH2}Let $(M,g)$ be a riemannian manifold. For $x\in
 M,\ R>0,$ we take a geodesic ball $B(x,R)$ such that, read in
 a chart $(V,\varphi ),\ B(x,R)\Subset V,$ the metric tensor
 at the center is the identity.\par 
We have a local Calderon Zygmund inequality on the manifold $M.$
 For any $r>1,$ there are constants $\displaystyle c_{1},c_{2}$
 depending only on $n=\mathrm{d}\mathrm{i}\mathrm{m}_{{\mathbb{R}}}M,\
 r$ and $R,$ such that:\par 
\quad \quad \quad \quad \quad $\forall u\in W^{2,r}(B(x,R)),\ {\left\Vert{u}\right\Vert}_{W^{2,r}(B(x,R/2))}\leq
 c_{1}{\left\Vert{u}\right\Vert}_{L^{r}(B(x,R))}+c_{2}{\left\Vert{\Delta
 u}\right\Vert}_{L^{r}(B(x,R))}.$
\end{lem}

\begin{proof}
We transcribe the problem in ${\mathbb{R}}^{n}$ by use of a coordinates
 path $(V,\varphi )$ exactly the same way we did to prove Theorem~\ref{RSMA0}.
 The Hodge laplacian is the second order elliptic matrix operator
 $\Delta _{\varphi }$ with ${\mathcal{C}}^{\infty }$ coefficients
 operating in $\varphi (V)\subset {\mathbb{R}}^{n}.$ By the choice
 of a $R$ small enough we can apply Lemma~\ref{RSMH4}, to the
 euclidean balls $B':=B_{e}(0,R_{e}')\subset \varphi (B(x,R/2)),\
 B:=B_{e}(0,R_{e})\subset \varphi (B(x,R))$ and we get, with
 $u_{\varphi }$ the $p$-form $u$ read in the chart $(V,\varphi ),$\par 
\[ {\left\Vert{u_{\varphi }}\right\Vert}_{W^{2,r}(B')}\leq c_{1}R^{-2}{\left\Vert{u_{\varphi
 }}\right\Vert}_{L^{r}(B)}+c_{2}{\left\Vert{\Delta _{\varphi
 }u_{\varphi }}\right\Vert}_{L^{r}(B)}.\]\par 
\quad The Lebesgue measure on $U$ and the canonical measure $dv_{g}$
 on $B(x,R)$ are equivalent; so the Lebesgue estimates and the
 Sobolev estimates up to order 2 on $U$ are valid in $B(x,R)$
 up to a constant.\par 
\quad So passing back to $M,$ we get, with $A:=\varphi ^{-1}(B),\ A':=\varphi
 ^{-1}(B')$\par 
\[ {\left\Vert{u}\right\Vert}_{W^{2,r}(A')}\leq c_{1}{\left\Vert{u}\right\Vert}_{L^{r}(A)}+c_{2}{\left\Vert{\Delta
 u}\right\Vert}_{L^{r}(A)}).\]\par 
So taking a smaller ball centered at $x$ we end the proof of the lemma.
\end{proof}

\begin{thm}
Let $M$ be a compact ${\mathcal{C}}^{\infty }$ Riemannian manifold
 without boundary. We have:\par 
\quad \quad \quad $\displaystyle \forall \omega \in L^{r}_{p}(M)\cap {\mathcal{H}}_{p}(M)^{\perp
 },\ r\in (1,\frac{n}{2});\ \exists u\in L^{s}_{p}(M)::\Delta u=\omega ,$\par 
with $\displaystyle \frac{1}{s}=\frac{1}{r}-\frac{2}{n}.$ Moreover
 $u\in W^{2,r}(M),\ {\left\Vert{u}\right\Vert}_{W^{2,r}(M)}\leq
 c{\left\Vert{\omega }\right\Vert}_{L^{r}(M)}.$
\end{thm}

\begin{proof}
First by duality we get the range $r>2.$ For this we shall proceed
 as we did in~\cite{AmarSt13}, using an avatar of the Serre duality~\cite{Serre55}.
 We take $t$ as in corollary~\ref{CC3}.\par 
\quad Let $g\in L^{t'}_{p}(M)\cap {\mathcal{H}}_{p}(M)^{\perp },$ we
 want to solve $\Delta v=g,$ with $t'>2$ and $t'$ conjugate to $t.$\par 
We know by the previous part that, with $r\leq 2,$\par 
\quad \quad \quad \begin{equation}  \forall \omega \in L^{r}_{p}(M)\cap {\mathcal{H}}_{p}(M)^{\perp
 },\ \exists u\in L^{t}_{p}(M),\ \Delta u=\omega ,\ {\left\Vert{u}\right\Vert}_{L^{t}(M)}\leq
 c{\left\Vert{\omega }\right\Vert}_{L^{r}(M)}.\label{HD5}\end{equation}\par 
Consider the linear form\par 
\quad \quad \quad $\forall \omega \in L^{r}_{p}(M),\ {\mathcal{L}}(\omega ):={\left\langle{u,g}\right\rangle},$\par
 
where $u$ is a solution of~(\ref{HD5}); in order for ${\mathcal{L}}(\omega
 )$ to be well defined, we need that if $u'$ is another solution
 of $\Delta u'=\omega ,$ then ${\left\langle{u-u',g}\right\rangle}=0;$
 hence we need that $g$ must be "orthogonal" to $p$-forms $\varphi
 $ such that $\Delta \varphi =0$ which is precisely our assumption.\par 
\quad Hence we have that ${\mathcal{L}}(f)$ is well defined and linear; moreover\par 
\quad \quad \quad $\ \left\vert{{\mathcal{L}}(f)}\right\vert \leq {\left\Vert{u}\right\Vert}_{L^{t}(M)}{\left\Vert{g}\right\Vert}_{L^{t'}(M)}\leq
 c{\left\Vert{\omega }\right\Vert}_{L^{r}(M)}{\left\Vert{g}\right\Vert}_{L^{t'}(M)}.$\par
 
So this linear form is continuous on $\omega \in L^{r}_{p}(M)\cap
 {\mathcal{H}}_{p}(M)^{\perp }.$ By the Hahn Banach Theorem there
 is a form $v\in L^{r'}_{p}(M)$ such that:\par 
\quad \quad \quad $\forall \omega \in L^{r}_{p}(M)\cap {\mathcal{H}}_{p}(M)^{\perp
 },\ {\mathcal{L}}(\omega )={\left\langle{\omega ,v}\right\rangle}={\left\langle{u,g}\right\rangle}.$\par
 
But $\omega =\Delta u,$ so we have, because $\Delta $ is essentially
 self adjoint and $M$ is compact, ${\left\langle{\omega ,v}\right\rangle}={\left\langle{\Delta
 u,v}\right\rangle}={\left\langle{u,\Delta v}\right\rangle}={\left\langle{u,g}\right\rangle},$
 for any $u\in L^{t}_{p}(M)::\Delta u\in L^{r}_{p}(M).$ In particular
 for $u\in {\mathcal{C}}^{\infty }_{p}(M).$ Now the hypothesis
 (iii) gives that $\displaystyle \Delta v=g$ in $\displaystyle
 L^{t'}_{p}(M),$ with $v\in L^{r'}_{p}(M).$ So we get:\par 
\quad \quad \quad $\forall g\in L^{t'}_{p}(M)\cap {\mathcal{H}}_{p}(M)^{\perp },\
 \exists v\in L^{r'}_{p}(M),\ \Delta v=g.$\par 
\quad It remains to prove the moreover. The condition (i') is true
 by the local existence Theorem~\ref{RSMA0}. The condition (iv')
 is true in the case of the Hodge Laplacian on a compact boundary-less
 manifold by Morrey's results for $L^{2}(M),$ so we are done for $r\leq 2.$\par 
\quad For $r>2,$ Lemma~\ref{RSMH2} gives us by compactness that there
 is a smaller $R>0$ and bigger constants $c_{1},\ c_{2}$ such that:\par 
\quad \quad \quad \begin{equation}  \forall x\in M,\ {\left\Vert{u}\right\Vert}_{W^{2,r}(B(x,R/2))}\leq
 c_{1}{\left\Vert{u}\right\Vert}_{L^{r}(B(x,R))}+c_{2}{\left\Vert{\Delta
 u}\right\Vert}_{L^{r}(B(x,R))}.\label{RSMH3}\end{equation}\par 
We take for $u$ our global solution in $L^{r}_{p}(M).$ We have
 that $\Delta u=\omega \in L^{r}_{p}(M)$ hence we can apply the
 estimate~(\ref{RSMH3}) to $u$:\par 
\quad \quad \quad $\displaystyle \forall x\in M,\ {\left\Vert{u}\right\Vert}_{W^{2,r}(B(x,R/2))}\leq
 c_{1}{\left\Vert{u}\right\Vert}_{L^{r}(B(x,R))}+c_{2}{\left\Vert{\Delta
 u}\right\Vert}_{L^{r}(B(x,R))}\leq C{\left\Vert{\omega }\right\Vert}_{L^{r}(M)}.$\par
 
Now it remains to cover $M$ with a finite set of balls of the
 type $B(x,R/2)$ to end the proof.
\end{proof}

\subsection{The $L^{r}$ Hodge decomposition.}
\quad In order to deduce the Hodge decomposition from the existence
 of a good solution to the Poisson equation, we shall need a
 little bit more material.\ \par 

\subsubsection{Basic facts.}
\quad Let ${\mathcal{H}}_{p}^{2}$ be the set of harmonic $p$-forms
 in $L^{2}(M),$ i.e. $p$-form $\omega $ such that $\Delta \omega
 =0,$ which is equivalent here to $d\omega =d^{*}\omega =0.$\ \par 
\quad The classical $L^{2}$ theory of Morrey~\cite{Morrey66} gives,
 on a compact manifold $M$ without boundary:\ \par 
\quad \quad \quad ${\mathcal{H}}_{p}:={\mathcal{H}}_{p}^{2}\subset {\mathcal{C}}^{\infty
 }(M)$ [~\cite{Morrey66}, (vi) p. 296]\ \par 
\quad \quad \quad $\mathrm{d}\mathrm{i}\mathrm{m}_{{\mathbb{R}}}{\mathcal{H}}_{p}<\infty
 $ [~\cite{Morrey66}, Theorem7.3.1].\ \par 
\quad This gives the existence of a linear projection from $L^{r}_{p}(M)\rightarrow
 {\mathcal{H}}_{p}$:\ \par 
\quad \quad \quad $\displaystyle \forall v\in L^{r}_{p}(M),\ H(v):=\sum_{j=1}^{N}{{\left\langle{v,e_{j}}\right\rangle}e_{j}}$\
 \par 
where $\lbrace e_{j}\rbrace _{j=1,...,N}$ is an orthonormal basis
 for ${\mathcal{H}}_{p}.$ This is meaningful because $v\in L^{r}_{p}(M)$
 can be integrated against $e_{j}\in {\mathcal{H}}_{p}\subset
 {\mathcal{C}}^{\infty }(M).$ Moreover we have $v-H(v)\in L^{r}_{p}(M)\cap
 {\mathcal{H}}_{p}^{\perp }$ in the sense that $\forall h\in
 {\mathcal{H}}_{p},\ {\left\langle{v-H(v),\ h}\right\rangle}=0;$
 it suffices to test on $h:=e_{k}.$ We get\ \par 
\quad \quad \quad $\displaystyle {\left\langle{v-H(v),\ e_{k}}\right\rangle}={\left\langle{v,e_{k}}\right\rangle}-{\left\langle{\sum_{j=1}^{N}{{\left\langle{v,e_{j}}\right\rangle}e_{j},e_{k}}}\right\rangle}={\left\langle{v,e_{k}}\right\rangle}-{\left\langle{v,e_{k}}\right\rangle}=0.$\
 \par 
\quad Let $v\in L^{r}_{p}(M).$ Set $h:=H(v)\in {\mathcal{H}}_{p},$
 and $\omega :=v-h.$ We have that $\forall k\in {\mathcal{H}}_{p},\
 {\left\langle{\omega ,k}\right\rangle}={\left\langle{v-H(v),\
 k}\right\rangle}=0.$ Hence we can solve $\Delta u=\omega $ with
 $u\in W^{2,r}_{p}(M)\cap L^{s}_{p}(M).$ So we get $v=h+\Delta
 u$ which means:\ \par 
\quad \quad \quad $\displaystyle L^{r}_{p}(M)={\mathcal{H}}_{p}^{r}\oplus \mathrm{I}\mathrm{m}\Delta
 (W^{2,r}_{p}(M)).$\ \par 
We have a direct decomposition because if $\omega \in {\mathcal{H}}_{p}\cap
 \mathrm{I}\mathrm{m}\Delta (W^{2,r}_{p}(M)),$ then $\omega \in
 {\mathcal{C}}^{\infty }(M)$ and\ \par 
\quad \quad \quad $\omega =\Delta u\Rightarrow \forall k\in {\mathcal{H}}_{p},\
 {\left\langle{\omega ,k}\right\rangle}=0$\ \par 
so choosing $k=\omega \in {\mathcal{H}}_{p}$ we get $\displaystyle
 {\left\langle{\omega ,\omega }\right\rangle}=0$ hence $\omega =0.$\ \par 
\quad Now we are in position to prove:\ \par 

\begin{thm}
If $(M,g)$ is a compact riemannian manifold without boundary;
 we have the strong $\displaystyle L^{r}$ Hodge decomposition:\par 
\quad \quad \quad $\displaystyle \forall r,\ 1\leq r<\infty ,\ L^{r}_{p}(M)={\mathcal{H}}_{p}^{r}\oplus
 \mathrm{I}\mathrm{m}d(W^{1,r}_{p}(M))\oplus \mathrm{I}\mathrm{m}d^{*}(W^{1,r}_{p}(M)).$
\end{thm}

\begin{proof}

 We already have\par 
\quad \quad \quad $L^{r}_{p}(M)={\mathcal{H}}_{p}^{r}\oplus \mathrm{I}\mathrm{m}\Delta
 (W^{2,r}_{p}(M))$\par 
where $\oplus $ means uniqueness of the decomposition.\par 
So: $\forall \omega \in L^{r}_{p}(M)\cap {\mathcal{H}}_{p}^{r\perp
 },\ \exists u\in W^{2,r}_{p}(M)::\Delta u=\omega .$\par 
From $\Delta =dd^{*}+d^{*}d$ we get $\omega =d(d^{*}u)+d^{*}(du)$
 with $du\in W^{1,r}_{p}(M)$ and $d^{*}u\in W^{1,r}_{p}(M),$ so\par 
\quad \quad \quad $\displaystyle \forall \omega \in L^{r}_{p}(M)\cap {\mathcal{H}}_{p}^{r\perp
 },\ \exists \alpha \in W^{1,r}_{p-1}(M),\ \exists \beta \in
 W^{1,r}_{p+1}(M)::\omega =\alpha +\beta ,$\par 
simply setting $\alpha =d(d^{*}u),\ \beta =d^{*}(du),$ which
 gives $\alpha \in d(W^{1,r}_{p-1}(M))$ and $\beta \in d^{*}(W^{1,r}_{p+1}(M))$
 and this proves the existence.\par 
\quad The uniqueness is given by Lemma 6.3 in~\cite{Scott95} and I
 copy this simple (but nice) proof here for the reader's convenience.\par 
\quad Suppose that $\alpha \in W^{1,r}_{p-1}(M),\ \beta \in W^{1,r}_{p+1}(M),\
 h\in {\mathcal{H}}_{p}$ satisfy $d\alpha +d^{*}\beta +h=0.$\par 
Let $\varphi \in {\mathcal{C}}^{\infty }_{p}(M),$ because of
 the classical ${\mathcal{C}}^{\infty }$-Hodge decomposition,
 there are $\eta \in {\mathcal{C}}^{\infty }_{p-1}(M),\ \omega
 \in {\mathcal{C}}^{\infty }_{p+1}(M)$ and $\tau \in {\mathcal{H}}_{p}$
 satisfying $\varphi =d\eta +d^{*}\omega +\tau .$\par 
Notice that $\ {\left\langle{d^{*}\beta ,d\eta }\right\rangle}={\left\langle{\beta
 ,d^{2}\eta }\right\rangle}={\left\langle{\beta ,0}\right\rangle}=0$
 and $\ {\left\langle{h,d\eta }\right\rangle}={\left\langle{d^{*}h,\eta
 }\right\rangle}=0,$ by the duality between $d$ and $d^{*}.$
 Linearity then gives\par 
\quad \quad \quad \begin{equation} {\left\langle{ d\alpha ,d\eta }\right\rangle}={\left\langle{d\alpha
 +d^{*}\beta +h,d\eta }\right\rangle}={\left\langle{0,d\eta }\right\rangle}=0.\label{CC0}\end{equation}\par
 
Finally we have\par 
\quad \quad \quad $\displaystyle {\left\langle{d\alpha ,\varphi }\right\rangle}={\left\langle{d\alpha
 ,d\eta }\right\rangle}+{\left\langle{d\alpha ,d^{*}\omega }\right\rangle}+{\left\langle{d\alpha
 ,\tau }\right\rangle}$\par 
\quad \quad \quad \quad \quad \quad \quad $\displaystyle =0+{\left\langle{\alpha ,d^{*2}\omega }\right\rangle}+{\left\langle{\alpha
 ,d\tau }\right\rangle}$    by~(\ref{CC0})\par 
\quad \quad \quad \quad \quad \quad \quad $\displaystyle ={\left\langle{\alpha ,0}\right\rangle}+{\left\langle{\alpha
 ,0}\right\rangle}$      because $\displaystyle d^{*2}=0$ and
 $\tau \in {\mathcal{H}}_{p}$\par 
\quad \quad \quad \quad \quad \quad \quad $\displaystyle =0.$\par 
Since ${\mathcal{C}}^{\infty }_{p}(M)$ is dense in $L^{r'}_{p}(M),\
 r'$ being the conjugate exponent of $r,$ and $\varphi $ is arbitrary,
 we see that $d\alpha =0.$ Analogously, we see that $d^{*}\beta
 =0$ and it follows that $h=0.$
\end{proof}

\section{Case of $\Omega $ a domain in $M.$}
\quad Let $\Omega $ be a domain in a ${\mathcal{C}}^{\infty }$ smooth
 complete riemannian manifold $M,$ compact or non compact; we
 want to show how the results in case of a compact boundary-less
 manifold apply to this case.\ \par 
\quad A classical way to get rid of a "annoying boundary" of a manifold
 is to use its "double". For instance: Duff~\cite{Duff52}, H\"ormander~\cite[p.
 257]{Hormand94}. Here we copy the following construction from~\cite[Appendix
 B]{GuneysuPigola}.\ \par 
\quad Let $N$ be a relatively compact domain of $M$ such that $\partial
 N$ is a smooth hypersurface and $\bar \Omega \subset N.$ The
 "Riemannian double" $D:=D(N)$ of $N,$ obtained by gluing two
 copies of $N$ along $\partial N,$ is a compact Riemannian manifold
 without boundary. Moreover, by its very construction, it is
 always possible to assume that $\displaystyle D$ contains an
 isometric copy of the original domain $N,$ hence of the original
 $\Omega .$ We shall also write $\Omega $ for its isometric copy
 to ease notations.\ \par 
\ \par 
\quad We shall need the following difficult result  by N. Aronszajn,
 A. Krzywicki and J. Szarski~\cite{Aronszajn62}, a strong continuation
 property, which says that if, for any compact set $K\subset
 M,$ the $p$-form $\omega $ satisfies the following inequality\ \par 
\quad \quad \quad \begin{equation} {\left\langle{ d\omega ,d\omega }\right\rangle}+{\left\langle{d^{*}\omega
 ,d^{*}\omega }\right\rangle}\leq C(K){\left\langle{\omega ,\omega
 }\right\rangle}\label{RSMH0}\end{equation}\ \par 
uniformly on $K,$ then if $\omega $ is zero to infinite order
 at a point $x_{0}\in M,$ we have that $\omega \equiv 0.$ The
 regularity conditions on $\omega $ are to be $L^{2}(M)$ with
 strong $L^{2}$ derivatives. The $p$-form $\omega $ must vanish
 at $x_{0}$ with all derivatives in "$L^{1}$ mean", which is
 also much weaker than the usual notion.\ \par 
We shall apply it to the compact manifold $D$ and with a ${\mathcal{C}}^{\infty
 }$ harmonic $p$-form $h,$ hence which satisfies inequality~(\ref{RSMH0}),
 and which is zero on an open non void set, which also implies
 a zero of infinite order.\ \par 
\quad The main lemma of this section is:\ \par 

\begin{lem}
~\label{tD0}Let $\omega \in L^{r}_{p}(\Omega ),$ then we can
 extend it to $\omega '\in L^{r}_{p}(D)$ such that: $\forall
 h\in {\mathcal{H}}_{p}(D),\ {\left\langle{\omega ',h}\right\rangle}_{D}=0.$
\end{lem}

\begin{proof}
Recall that ${\mathcal{H}}_{p}(D)$ is the vector space of harmonic
 $p$-form in $D,$ it is of finite dimension $K_{p}$ and ${\mathcal{H}}_{p}(D)\subset
 {\mathcal{C}}^{\infty }(D).$\par 
\quad Make an orthonormal basis $\lbrace e_{1},...,e_{K_{p}}\rbrace
 $ of ${\mathcal{H}}_{p}(D)$ with respect to $L^{2}_{p}(D),$
 by the Gram-Schmidt procedure.We get $\ {\left\langle{e_{j},e_{k}}\right\rangle}_{D}:=\int_{D}{e_{j}e_{k}dV}=\delta
 _{jk}.$\par 
Set $\lambda _{j}:={\left\langle{\omega {\11}_{\Omega },e_{j}}\right\rangle}={\left\langle{\omega
 ,e_{j}{\11}_{\Omega }}\right\rangle},\ j=1,...,K_{p},$ which
 makes sense since $e_{j}\in {\mathcal{C}}^{\infty }(D)\Rightarrow
 e_{j}\in L^{\infty }(D),$ because $D$ is compact.\par 
\par 
\quad We shall see that the system $\lbrace e_{k}{\11}_{D\backslash
 \Omega }\rbrace _{k=1,...,K_{p}}$ is a free one. Suppose this
 is not the case, then it will exist $\gamma _{1},...,\gamma
 _{K_{p}},$ not all zero, such that $\sum_{k=1}^{K_{p}}{\gamma
 _{k}e_{k}{\11}_{D\backslash \Omega }}=0$ in $D\backslash \Omega
 .$ But the function $\displaystyle h:=\sum_{k=1}^{K_{p}}{\gamma
 _{k}e_{k}}$ is in ${\mathcal{H}}_{p}(D)$ so if $h$ is zero in
 $D\backslash \Omega $ which is non void, then $h\equiv 0$ in
 $D$ by the N. Aronszajn, A. Krzywicki and J. Szarski~\cite{Aronszajn62}
 result. This is not possible because the $e_{k}$ make a basis
 for ${\mathcal{H}}_{p}(D).$ So the system $\lbrace e_{k}{\11}_{D\backslash
 \Omega }\rbrace _{k=1,...,K_{p}}$ is a free one.\par 
\quad We set $\gamma _{jk}:={\left\langle{e_{k}{\11}_{D\backslash \Omega
 },e_{j}{\11}_{D\backslash \Omega }}\right\rangle}$ hence we
 have that $\mathrm{d}\mathrm{e}\mathrm{t}\lbrace \gamma _{jk}\rbrace
 \neq 0.$ So we can solve the linear system to get $\lbrace \mu
 _{k}\rbrace $ such that\par 
\quad \quad \quad \begin{equation}  \forall j=1,...,K_{p},\ \sum_{k=1}^{K_{p}}{\mu
 _{k}{\left\langle{e_{k}{\11}_{D\backslash \Omega },e_{j}}\right\rangle}}=\lambda
 _{j}.\label{RSMH1}\end{equation}\par 
We put $\omega '':=\sum_{j=1}^{K_{p}}{\mu _{j}e_{j}{\11}_{D\backslash
 \Omega }}$ and $\displaystyle \omega ':=\omega {\11}_{\Omega
 }-\omega ''{\11}_{D\backslash \Omega }=\omega -\omega ''.$ From~(\ref{RSMH1})
 we get\par 
\quad \quad \quad $\displaystyle \forall j=1,...,K_{p},\ {\left\langle{\omega ',e_{j}}\right\rangle}_{D}={\left\langle{\omega
 ,e_{j}}\right\rangle}-{\left\langle{\omega '',e_{j}}\right\rangle}=\lambda
 _{j}-\sum_{k=1}^{K_{p}}{\mu _{k}{\left\langle{e_{k}{\11}_{D\backslash
 \Omega },e_{j}}\right\rangle}}=0.$\par 
So the $p$-form $\omega '$ is orthogonal to ${\mathcal{H}}_{p}.$
 Moreover $\omega '_{\mid \Omega }=\omega $ and clearly $\omega
 ''\in L^{r}_{p}(D)$ being a finite combination of $\displaystyle
 e_{j}{\11}_{D\backslash \Omega },$ so $\omega '\in L^{r}_{p}(D)$
 because $\omega $ itself is in $\displaystyle L^{r}_{p}(D).$
 The proof is complete.
\end{proof}
\ \par 
\quad Now let $\omega \in L^{r}(\Omega )$ and see $\Omega $ as a subset
 of $D;$ then extend $\omega $ as $\omega '$ to $D$ by Lemma~\ref{tD0}.\ \par 
By the results on the compact manifold $D,$ because $\omega '\perp
 {\mathcal{H}}_{p}(D),$ we get that there exists $u'\in W^{2,r}_{p}(D),\
 u'\perp {\mathcal{H}}_{p}(D),$ such that $\Delta u'=\omega ';$
 hence if $u$ is the restriction of $u'$ to $\Omega $ we get
  $u\in W^{2,r}_{p}(\Omega ),\ \Delta u=\omega $ in $\Omega .$\ \par 
Hence we proved\ \par 

\begin{thm}
~\label{tD1}Let $\Omega $ be a domain in the smooth complete
 riemannian manifold $M$ and $\omega \in L^{r}_{p}(\Omega ),$
 then there is a $p$-form $u\in W^{2,r}_{p}(\Omega ),$ such that
 $\Delta u=\omega $ and ${\left\Vert{u}\right\Vert}_{W^{2,r}_{p}(\Omega
 )}\leq c(\Omega ){\left\Vert{\omega }\right\Vert}_{L^{r}_{p}(\Omega )}.$
\end{thm}
\quad As for domains in ${\mathbb{R}}^{n},$ there is no constrain for
 solving the Poisson equation in this case.\ \par 

\begin{cor}
(Of the proof) Let $M$ be a smooth compact riemannian manifold
 with smooth boundary $\partial M.$ Let $\displaystyle \omega
 \in L^{r}_{p}(M).$ There is a form $u\in W^{2,r}_{p}(M),$ such
 that $\Delta u=\omega $ and ${\left\Vert{u}\right\Vert}_{W^{2,r}_{p}(M)}\leq
 c{\left\Vert{\omega }\right\Vert}_{L^{r}_{p}(M)}.$
\end{cor}

\begin{proof}
We can build the "double manifold" $D:=D(M)$ which is compact
 without boundary. Copying the proof of Theorem~\ref{tD1} we
 extend the form $\omega $ defined on $M,$ viewed as a subset
 in $D,$ to a form $\omega '$ in $\displaystyle L^{r}_{p}(D)$
 orthogonal to ${\mathcal{H}}_{p}(D)$ so there is a $u'::\Delta
 u'=\omega '$ in $D.$ We just take $u:=u'_{\mid M}$ to finish the proof.
\end{proof}
\ \par 

\bibliographystyle{/usr/local/texlive/2013/texmf-dist/bibtex/bst/base/siam}

\end{document}